\documentclass[12pt]{article}
\usepackage{graphicx}					
\usepackage{mathrsfs}
\usepackage{verbatimbox}
\usepackage{listofitems}
\usepackage{amsmath, amssymb}
\usepackage{amsfonts, ifthen, amsthm, vmargin}
\usepackage[usenames]{color}
\usepackage{ulem}
\usepackage[colorlinks=true, urlcolor=blue2, citecolor=cite]{hyperref}

\DeclareFontFamily{OT1}{rsfs}{} \DeclareFontShape{OT1}{rsfs}{m}{n}{
<-7> rsfs5 <7-10> rsfs7 <10-> rsfs10}{}
\DeclareMathAlphabet\mathcurl{OT1}{rsfs}{m}{n}

\theoremstyle{plain}
\newtheorem*{cor*}{Corollary}

\theoremstyle{plain}
\newtheorem*{gt*}{Grothendieck's Inequality (Lindenstrauss-Pe{\l}czy\'{n}ski style)}

\theoremstyle{plain}
\newtheorem*{thmGI*}{Grothendieck's Identity}

\theoremstyle{plain}
\newtheorem*{thmHI*}{Haagerup's Identity}

\theoremstyle{definition}
\newtheorem*{df*}{Definition}

\theoremstyle{definition}
\newtheorem*{ex*}{Example}

\theoremstyle{definition}
\newtheorem*{ack*}{Acknowledgments}

\theoremstyle{plain}
\newtheorem*{exkrv*}{Example (Krivine's constant)}

\theoremstyle{plain}
\newtheorem*{gi*}{Matrix Version of Grothendieck's Identity}

\theoremstyle{plain}
\newtheorem*{guidex*}{Guiding Example}

\theoremstyle{plain}
\newtheorem*{prp*}{Proposition}

\theoremstyle{plain}
\newtheorem*{cprp*}{Proposition (correlation matrix version of GT and little GT)}

\theoremstyle{plain}
\newtheorem*{thm*}{Theorem}

\theoremstyle{plain}
\newtheorem*{thmgn*}{Theorem (A. Grothendieck, 1953; H. Niemi, 1983)}

\theoremstyle{plain}
\newtheorem*{c*}{Conjecture}

\theoremstyle{definition}
\newtheorem*{rem*}{Observation}

\theoremstyle{definition}
\newtheorem*{rem2*}{Remark}

\theoremstyle{plain}
\newtheorem*{arcsinex*}{Example (Stieltjes, 1889)}

\theoremstyle{plain}

\theoremstyle{plain}

\theoremstyle{plain}

\newcount\colveccount
\newcommand*\colvec[1]{
        \global\colveccount#1
        \begin{pmatrix}
        \colvecnext
}
\def\colvecnext#1{
        #1
        \global\advance\colveccount-1
        \ifnum\colveccount>0
                \\
                \expandafter\colvecnext
        \else
                \end{pmatrix}
        \fi
}

\definecolor{blue2}{rgb}{0.08, 0.38, 0.74}
\definecolor{cite}{rgb}{1.0, 0.13, 0.32}

\newcommand\lb{\left(}
\newcommand\rb{\right)}
\newcommand{\lsb}{\left[}
\newcommand{\rsb}{\right]}
\newcommand\ignore[1]{}
\newcommand\A{\mathbb A}
\newcommand\C{\mathbb C}
\newcommand\D{\mathbb D}
\newcommand\R{\mathbb R}
\newcommand\T{\mathbb T}
\newcommand\M{\mathbb M}
\newcommand\N{\mathbb N}
\newcommand{\F}{\mathbb{F}}
\newcommand\E[2][\P]{\mathbb E_{#1}\left[ #2\right]}
\renewcommand\E{\mathbb E}
\renewcommand\P{\mathbb P}

\renewcommand{\S}{\mathbb{S}}

\newcommand\ind{1\hspace{-2.5mm}1}

\newcommand\madd[1]{{\color{blue}#1}}

\newcommand\fadd[1]{{\color{red}#1}}

\makeatletter
\renewcommand{\@seccntformat}[1]{\csname the#1\endcsname.\hspace{1em}}%

\makeatother

\begin{document}

\title{Grothendieck's inequality and completely correlation preserving functions - a summary 
of recent results and an indication of related research problems
\thanks{\textit{MSC2020 subject classifications.} Primary 05A10, 33C05, 41A58, 62H05, 62H20; 
secondary 15A60, 46A20, 47B10.}}
\author{Frank Oertel  \\
	Philosophy, Logic \& Scientific Method\\
	Centre for Philosophy of Natural and Social Sciences (CPNSS)\\
	London School of Economics and Political Science\\
	Houghton Street, London WC2A 2AE, UK\\
	\texttt{f.oertel@email.de}
	}

\date{} 

\maketitle


\noindent {\bf{Abstract.}} 
As part of the search for the value of the smallest upper bound of the best constant for the 
famous Grothendieck inequality, the so-called Grothendieck constant (a hard open problem - unsolved since 1953), we provide a further approach, primarily built on functions which 
map correlation matrices entrywise to correlation matrices by means of the Schur product, 
multivariate Gaussian analysis, copulas and inversion of suitable Taylor series. We summarise 
first results and point towards related open problems and topics for future research.   
\section{Introduction}
Despite its emergence more than six decades ago, the techniques and results of the 
actually pathbreaking work of A. Grothendieck in the metric theory of tensor products are 
still not widely known nor appreciated. Very likely this is due to the fact that Grothendieck 
included virtually no proofs, and that he used the (duality) theory of the rather abstract 
(yet very powerful) notion of tensor products of Banach spaces (cf. \cite{G1953, DF1993, 
DFS2008}). The fundamental idea exploited in \cite{DF1993} is a one-to-one correspondence 
between Grothendieck's finitely generated tensor norms and maximal Banach operator ideals (in 
the sense of Pietsch - cf. \cite{P1980}) via trace duality. Theory and applications of 
operator ideals are widely known (not by functional analysts only), as opposed to the tensor 
norm theory of Grothendieck so that also \cite{DF1993} (such as \cite{DFS2008}) is a very 
valuable source which strongly helps to make Grothendieck's approach accessible to a wider 
community. 
\\[0.2em]
In particular, the famous Grothendieck inequality (also known as \textit{the fundamental theorem of the metric theory of tensor products}), published in Grothendieck's famous paper 
\cite{G1953} 
had a profound influence on the geometry of Banach spaces and operator theory in the 1970s and 
1980s. Meanwhile, in addition to this impact, Grothendieck's inequality exhibits deep  
applications in different directions (including theoretical computer science, 
computational complexity, analysis of Boolean functions, random graphs (including 
the mathematics of the systemic risk in financial networks, analysis of nearest-neighbour 
interactions in a crystal structure (Ising model), correlation clustering and image 
segmentation in the field of computer vision), NP-hard combinatorial optimisation, 
non-convex optimisation and semidefinite programming (cf. \cite{GW1995}), foundations and 
philosophy of quantum mechanics, quantum information theory, quantum correlations, quantum 
cryptography, communication complexity protocols and even high-dimensional private data 
analysis (cf. \cite{DNT2015})! Also in these fields it offers many challenging related open 
questions.
\\[0.2em]
The interest in Grothendieck's work revived when J. Lindenstrauss and A. 
Pe{\l}czy\'{n}ski recast its main results in the more traditional language of operators and 
matrices (see \cite{LP1968} and \cite[Theorem A.3.1]{DFS2008}) which is also the basis 
of our own research. A slightly bit modified version of this rewritten version of 
Grothendieck's inequality reads as follows (cf. \cite[Lemma 2.2]{FLZ2018}):

\begin{gt*}
Let $\F \in \{\R, \C\}$. There is an absolute constant $K>0$ such that for 
any $m, n \in \N$, for any $m \times n$ matrix $\big(a_{ij}\big)$ with entries in 
$\F$, any $\F$-Hilbert space $H$ and any vectors $u_1, \ldots, u_m, v_1, \ldots, v_n \in H$ 
of norm $=1$ the following inequality is satisfied:
\[
\Big\vert \sum_{i = 1}^{m}\sum_{j = 1}^{n} a_{ij}\langle u_i, v_j\rangle_H \Big\vert  \leq K 
\, \sup{\Big\{ \Big\vert \sum_{i = 1}^{m}\sum_{j = 1}^{n} a_{ij} p_i q_j\Big\vert : 
\vert p_i \vert = 1 = \vert q_j \vert \, \forall i,j\Big\}}\,.
\]
\end{gt*}
\noindent The smallest possible value of the corresponding constant $K$ is called the 
\textit{Grothendieck constant $K_G^\F$}. The superscripts $\R$ and $\C$ are used to indicate 
the different values in the real and complex case. The complex constant is smaller than the 
real one, since (cf. \cite{K1977, P2012}):
\begin{thm*}
\[
K_G^\C \leq e^{1-\gamma} < \frac{\pi}{2} < K_G^\R \leq \sqrt{2}\,K_G^\C\,,
\]
where $\gamma : = \sum_{n=2}^\infty (-1)^n \frac{\zeta(n)}{n} = -\Gamma'(1) 
\approx 0,577215664901533 \ldots$ denotes the Euler-Mascheroni constant.
\end{thm*}
\noindent Computing the exact numerical value of the constants $K_G^\R$ and $K_G^\C$ is still 
an open problem (unsolved since 1953). Here, our own research activities are attached. We 
look for a general framework (primarily build on methods originating from (block) matrix 
analysis (cf. \cite{HJ2013}), multivariate statistics with real and complex Gaussian 
random vectors, theory of special functions, modelling of statistical dependence with 
copulas and combinatorics, whose complexity increases rapidly in dimension, though) which 
allows us either to give the value of $K_G^\R$, respectively $K_G^\C$ explicitly or to 
approximate 
these values from above and from below at least. Surprisingly, our approach - which in 
particular allows a short proof of the real Grothendieck inequality, even with Krivine's 
upper bound of $K_G^\R$ - confronts us strongly with the question whether the seemingly non-
avoidable combinatoric complexity actually allows us to determine the values of $K_G^\R$, 
respectively $K_G^\F$ explicitly, or not (cf. research problem 7.1 below).   
\\[0.2em]
As usual, $\N_0 : = \N \cup \{0\}$ denotes the set of all non-negative integers. Additionally, for any $k \in \N$, we put $\N_k : = \{n : n \in \N \text{ and } n \geq k\}$. In 
the case of finite-dimensional Hilbert spaces Grothendieck's inequality in matrix form is 
given by the following well-known result: 
\begin{thm*}
Let $\F \in \{\R, \C\}$. Then for any $d \in \N$ there is a constant $K(d) > 0$ such 
that
\[
\Big\vert \sum_{i = 1}^{m}\sum_{j = 1}^{n} a_{ij}\langle u_i, v_j\rangle_{l_2^d} 
\Big\vert  \leq K(d) \, \sup{\Big\{ \Big\vert \sum_{i = 1}^{m}\sum_{j = 1}^{n} a_{ij} 
p_i q_j\Big\vert : \vert p_i \vert = 1 = \vert q_j \vert \, \forall i,j \Big\}}\,,
\] 
for any $m, n \in \N$, for any matrix $\lb a_{ij} \rb \in \M(m \times n;\F)$, for any 
vectors $u_1, \ldots, u_m, v_1, \ldots, v_n \in \S^{d-1}$. Here, $\S^{d-1} : = \{w : 
w \in l_2^d \text{ and } \Vert w \Vert_2 = 1\}$ denotes the unit sphere of $l_2^d \equiv 
l_2^d(\F)$.
\end{thm*}
\begin{rem2*}
Note that we may include the case $d=1$ here, since if $d=1$, $\Vert w \Vert 
= \vert w \vert$ and $u^\ast v = v\overline{u}$ for all $u,v,w \in \F$, implying 
that by our definition $\S^{0} = \{z \in \F : \vert z \vert = 1\}$.
\end{rem2*}
\noindent Let $K_G^\F(d)$ denote the smallest possible value 
of the corresponding constant $K(d)$. 
Since the sequence $\big(K_G^\F(d) \big)_{d \in \N}$ obviously is non-decreasing it follows 
that $K_G^\F = \lim\limits_{d \to \infty}K_G^\F(d) = \sup\{K_G^\F(d): d \in \N\}$. 
Moreover, we may add: 
\begin{prp*}
Let $d \in \N$. Then
\[
K_G^\R(2d) \leq \sqrt{2}\,K_G^\C(d)\,.
\]
\end{prp*}
\noindent In particular, by taking the limit $d \to \infty$, we reobtain $K_G^\R \leq \sqrt{2}
\,K_G^\C$. 
\\[0.2em]
An important special case of Grothendieck's inequality (known as \textit{little Grothendieck 
inequality}) appears if just positive semidefinite matrices $A$ are considered. Let $k_G^\F$ 
denote the Grothendieck constant, derived from Grothendieck's inequality restricted to the set 
of all positive semidefinite $n \times n$ matrices over $\F$. Then (cf. \cite{N1983}):
\begin{thm*}[A. Grothendieck, 1953; H. Niemi, 1983]
Let $\F \in \{\R, \C\}$ and $H$ be an arbitrary Hilbert space over $\F$. Let $n \in \N$. Then
\begin{itemize}
\item[$($i$)$] 
\[
k_G^\R = \frac{\pi}{2} \hspace{0.5cm}\text{ $($A. Grothendieck$)$ }
\] 
\end{itemize}
and
\begin{itemize}
\item[$($ii$)$] 
\[
k_G^\C = \frac{4}{\pi} \hspace{0.5cm}\text{ $($H. Niemi$)$ }.
\]
\end{itemize}
In particular, $\frac{\pi}{2} \leq K_G^\R$ and $\frac{4}{\pi} \leq K_G^\C$.
\end{thm*}
\noindent Until present the following encapsulation of $K_G^\R$ holds (cf. \cite{BMMN2013} 
and the related references therein (including \cite{K1977})):
\[
1,676 < K_G^\R \,\, \stackrel{(!)}{<} \,\, \frac{\pi}{2 \ln(1 + \sqrt{2})} \approx 1,782\,.
\]
\section{A correlation matrix version of Grothendieck's inequality}
Regarding our approach, the following equivalent reformulation of the Grothendieck inequality 
(in the following abbreviated by ``GT'') which discloses a link to correlation matrices 
(and hence to multivariate statistics and Gaussian analysis) is of crucial importance. 
Recall that for fixed $k \in \N$ a $k \times k$ correlation matrix precisely coincides with a 
positive semidefinite matrix $\Sigma \equiv \big(\sigma_{ij}\big)_{i, j \in [k]} \in 
\M(k \times k; \F)$ which in addition satisfies $\sigma_{ii} =1$ for all $i \in [k]$ 
(i.\,e., the diagonal of $\Sigma$ is filled with 1's only). Of particular relevance is the set 
$C_1(k; \F)$ of all $k \times k$ correlation matrices of rank 1 and the (Hermitian) block 
matrix $J(A) : = \frac{1}{2}\begin{pmatrix}
0 & A\\
A^\ast & 0
\end{pmatrix}$, where $m, n \in \N$ and $A \in M(m\times n; \F)$. Observe that in the 
following equivalent reformulation of GT - seemingly - no Hilbert space $H$ is needed\,!
\begin{cprp*}
Let $\F \in \{\R, \C\}$ and $m, n \in \N$. Let $\Sigma$ be an arbitrary $(m+n) \times (m+n)$ 
correlation matrix with entries in $\F$ and $A \in \M(m \times n; \F)$. Viewing $A$ as bounded 
linear operator from $l_\infty^n$ into $l_1^m$, we have
\[
\vert{\text{tr}}\big(J(A)\,\Sigma\big)\vert 
\leq K_G^\F \max_{\Theta \in C_1(m+n; \F)}\vert{\text{tr}}\big(J(A)\,\Theta\big)\vert = 
K_G^\F \Vert A \Vert_{\infty, 1}\,. 
\]
If in addition $A$ is positive semidefinite then the little GT is equivalent to 
\[
\vert{\text{tr}}\big(A^\ast\,\Sigma\big)\vert 
\leq k_G^\F \max_{\Theta \in C_1(m+n; \F)}\vert{\text{tr}}\big(A^\ast\,\Theta\big)\vert\,. 
\]
\end{cprp*}
\begin{rem*}
Let $k \in \N$. As we just have seen the \textit{little} GT can be equivalently written as
\[
\vert{\text{tr}}\big(B^\ast\,\Sigma\big)\vert 
\leq k_G^\F \max_{\Theta \in C_1(k; \F)}\vert{\text{tr}}\big(B^\ast\,\Theta\big)\vert = 
k_G^\F \Vert B \Vert_{\infty, 1}\,,
\] 
where $B \in \M(k \times k; \F)$ is an arbitrary \textit{positive semidefinite} matrix 
and $\Sigma$ an arbitrary $k \times k$ correlation matrix with entries in $\F$. 
However, we don't know whether we also may substitute in GT itself the block 
matrix $J(A)$ through an \textit{arbitrary} matrix $B \in \M(k \times k; \F)$. If this were 
the case an application of the bipolar theorem shows that the latter would be equivalent to 
the non-negligible set inclusion
\[
C(k; \F) \subseteq K_G^\F \, {}{acx\big(C_1(k; \F)\big)} 
\text{ for all } k \in \N\,,
\]
where $C(k; {}{\F})$ denotes the set of \textit{all} $k \times k$ correlation matrices with entries in $\F$ and $acx$ the absolute convex hull (of a set). It can be shown that
\[
C(k; \F) \subseteq K_G^\F \, {}{acx}\big(\{xy^\top : (x, y) \in 
(\S^{0})^k \times (\S^{0})^k\}\big) = K_G^\F \, B_{{}{\mathcal{N}(l_1^k, l_\infty^k)}} 
\text{ for all } k \in \N
\]
though, where $B_{{}{\mathcal{N}(l_1^k, l_\infty^k)}}$ denotes the unit ball of the Banach 
space of nuclear operators between $l_1^k$ and $l_\infty^k$, equipped with the nuclear norm
(which should not be mixed up with $B_{{}{\mathcal{N}(l_2^k, l_2^k)}} = 
cx\big(\{xy^\top : \Vert x\Vert_2\,\Vert y \Vert_2 = 1\}\big)$ (cf. {\url{https://
convexoptimization.com/TOOLS/0976401304.pdf}}, Example 2.3.2.0.2, equation (97))).
This non-trivial functional analytic result (whose proof involves the 
structure of the set of all \textit{quantum} correlation matrices (described by B. S. 
Tsirelson - cf. \cite{P2018} and cited references therein) and a description of the extreme 
points of $B_{{}{\mathcal{N}(l_1^n, l_\infty^m)}}$) is not subject of discussion 
in this document, though.
\end{rem*}
\noindent In fact, if we allow the implementation of a possibly strictly larger absolute 
constant than $K_G^\F$ it is possible to deduce a further non-trivial inequality - which even implies GT as a corollary! Namely,
\begin{thm*}
Let $\F \in \{\R, \C\}$. Then there exists an absolute constant $K_\ast^\F > 1$ such that
\[
\vert{\text{tr}}\big(B^\ast\,\Sigma\big)\vert \leq 
K_\ast^\F \max_{\Theta \in C_1(k; \F)}\vert{\text{tr}}\big(B^\ast\,\Theta\big)\vert\,. 
\]
for any $k \in \N$, any $\Sigma \in C(k; \F)$ and any $B \in \M(k \times k; \F)$. Moreover, 
\[
C(k; \F) \subseteq K_\ast^\F \,\text{acx}\lb C_1(k; \F) \rb \text{ for all } k \in \N\,,
\]
$K_\ast^\R \in [K_G^\R, \sinh\big(\frac{\pi}{2}\big)]$ and $K_\ast^\C \in 
[K_G^\C, \frac{8}{\pi}-1]$.
\end{thm*}
\section{Grothendieck's identity and Haagerup's idendity: a common source}
\noindent The main ingredients of the proof of GT are two \textit{equalities}, 
namely Grothendieck's identity (if $\F = \R$ - cf. e.g. the proof of 
\cite[Prop. 4.4.2]{DFS2008}) and Haagerup's identity (if $\F = \C$ - see \cite{H1987}). 
Rewritten in terms of real Gaussian random vectors (if $\F = \R$) and complex proper Gaussian 
random vectors (if $\F = \C$) they imply the following two results, 
revealing a common underlying structure 
for both fields, $\R$ and $\C$. To this end, recall (e.g. from \cite{AAR1999}) that for any 
$a,b\in \C$, any $c \in \C\setminus\{-n : n \in \N_0\}$ and any $z \in 
\D = \{z \in \C : \vert z \vert < 1\}$ the well-defined power series
\[
{}_2 F_1(a, b, c; z) : = \frac{\Gamma(c)}{\Gamma(a)\Gamma(b)}\sum_{n=0}^{\infty} 
\frac{\Gamma(a+n)\Gamma(b+n)}{\Gamma(c+n)}\,\frac{z^n}{n!}
\]
denotes the \textit{Gaussian hypergeometric function}. If in addition $\Re(c) > \Re(a+b)$ then 
the series converges absolutely on $\partial\D = \{z \in \C : \vert z \vert = 1\}$ and 
satisfies ${}_2 F_1(a, b, c; 1) = \frac{\Gamma(c)\Gamma(c-a-b)}{\Gamma(c-a)\Gamma(c-b)}$ 
(Gauss' Summation Theorem).
\begin{thmGI*}
Let $n \in \N$, $u, v \in \R^n$ such that $\Vert u \Vert = 1$ and $\Vert v \Vert = 1$ and 
$X \sim N_n(0, I_n)$ be a standard-normally distributed real Gaussian random vector. Then
\begin{eqnarray*}
\int_{{\S^{n-1}}}{\text{sign}}(u^\top x){\text{sign}}(v^\top x)\,
\underline{\sigma}^{n-1}(dx) & = & 
\E[{\text{sign}}(u^\top X){\text{sign}}(v^\top X)] \, = \, 
\frac{2}{\pi}\arcsin\lb u^\top v\rb\\
& = & \E[\vert X_1 \vert]^2\,u^\top v\,
{}_2 F_1\lb\frac{1}{2}, \frac{1}{2}, \frac{3}{2}; \vert u^\top v \vert^2 \rb\,,
\end{eqnarray*}
where 
${\underline{\sigma}}^{n-1} : = \frac{\Gamma(n/2)}{2 \pi^{n/2}}\,\sigma^{n-1}$ denotes 
the normalised (``uniform'') surface area probability measure on the unit sphere 
$\S^{n-1}$ and $\P(X \in \cdot) = \gamma_n$ the 
Gaussian probability measure on $\R^n$. 
\end{thmGI*}
\begin{rem2*}
If $n=1$, we put $\sigma^{0}(\S^{0}) = \sigma^{0}(\{-1,1\}): = 2$.
\end{rem2*} 
\begin{thmHI*}
Let $n \in \N$, $u, v \in \C^n$ such that $\Vert u \Vert = 1$ and $\Vert v \Vert = 1$. Let
$Z \sim {\C}N_n(0, I_n)$ be a standard-normally distributed complex proper 
Gaussian random vector. Then

\scalebox{0.9}{
\vbox{
\begin{eqnarray*}\label{eq:Haagerup}
\int_{{\C^n}}{\text{sign}}(u^\ast z){\text{sign}}(\overline{v^\ast z})\,\gamma_n^{(\C)}(dz) 
& = & \E[{\text{sign}}(u^\ast Z){\text{sign}}(\overline{v^\ast Z})] \, = \, 
{\text{sign}}(u^\ast v)\,
\frac{1}{4}\int_{0}^{2\pi}\arcsin\big(\vert u^\ast v \vert\cos(t)\big)\cos(t)\,dt\\
& = & \frac{\pi}{4}\,u^\ast v\,{}_2 F_1\lb\frac{1}{2}, 
\frac{1}{2}, 2; \vert u^\ast v \vert^2 \rb \, = \, 
\E[\vert Z_1 \vert]^2\,u^\ast v\,
{}_2 F_1\lb\frac{1}{2},\frac{1}{2}, 2; \vert u^\ast v \vert^2 \rb \,,
\end{eqnarray*}
}}

\noindent where $\gamma_n^{(\C)}(B) : = \gamma_{2n}^{(\R)}
\big(\big\{(\sqrt{2}\,\Re(z)^\top, \sqrt{2}\,\Im(z)^\top)^\top : z \in B \big\}\big) 
\hspace{0.1cm} (B \in {\mathcal{B}}(\C^n))$ denotes the Gaussian probability measure on 
$\C^n$.
\end{thmHI*}
\begin{rem*}
If in addition $u^\ast v \not=0$, we may add the following \textit{real} 
spherical integral representation in Haagerup's identity:

\scalebox{0.9}{
\vbox{
\[
\int_{{\C^n}}{\text{sign}}(u^\ast z){\text{sign}}(\overline{v^\ast z})\,\gamma_n^{(\C)}(dz)
= {\text{sign}}(u^\ast v)\Re\lb\E[f_{u,v}(Z)]\rb =
{\text{sign}}(u^\ast v)\,\int_{\S^{2n-1}} \Re\big(f_{u,v}(x+iy)\big)
\underline{\sigma}^{2n-1}(d(x,y))\,,
\]
}}

where $\C^n \ni z \mapsto f_{u,v}(z) : = 
\frac{1}{{\text{sign}}(u^\ast v)}\,{\text{sign}}(u^\ast z){\text{sign}}(\overline{v^\ast z})$. 
\end{rem*}
\noindent It can be shown that both identities arise as a special case of the following 
result, where we explicitly describe all non-negative integer powers of an expectation of 
inner products of suitably correlated - real - Gaussian random vectors. Here, we possibly 
should point to the so-called kernel trick, used also for the computation of inner products 
in high-dimensional feature spaces using simple functions defined on pairs of input patterns 
which is a crucial ingredient of support vector machines in statistical learning theory; i.e., 
learning machines that construct decision functions of sign type. This trick allows the 
formulation of nonlinear variants of any algorithm that can be cast in terms of inner 
products (cf. \cite[Chapter 5.6]{V2000}).
\\[0.2em]
To this end, let us consider the real correlation matrices
\[
\Sigma_{2d}(\rho) : = 
\begin{pmatrix}
    I_d & \rho\,I_d \\
    \rho\,I_d & I_d
  \end{pmatrix}
=
\begin{pmatrix}
   \begin{matrix}
    1 & 0 & \ldots & 0  \\
    0 & 1 & \ldots & 0 \\
		\vdots & \vdots & \ddots & \vdots \\
		0 & 0 & \ldots & 1
  \end{matrix}
					
		&

   \begin{matrix}
    \rho & 0 & \ldots & 0  \\
    0 & \rho & \ldots & 0 \\
		\vdots & \vdots & \ddots & \vdots \\
		0 & 0 & \ldots & \rho
  \end{matrix}
			
		\\

   \begin{matrix}
    \rho & 0 & \ldots & 0  \\
    0 & \rho & \ldots & 0 \\
		\vdots & \vdots & \ddots & \vdots \\
		0 & 0 & \ldots & \rho
  \end{matrix}

				& 
		
   \begin{matrix}
    1 & 0 & \ldots & 0  \\
    0 & 1 & \ldots & 0 \\
		\vdots & \vdots & \ddots & \vdots \\
		0 & 0 & \ldots & 1
  \end{matrix}
	
	\end{pmatrix}\,,
\]
where $-1 \leq \rho \leq 1$ and $d \in \N$. 
\begin{thm*}
Let $d, m \in \N$, $\rho \in (-1,1)$ and $(X^\top, Y^\top)^\top \equiv 
(X_1, \ldots, X_d, Y_1, \ldots, Y_d)^\top \sim N_{2d}\big(0, \Sigma_{2d}(\rho)\big)$. 
\begin{itemize}
\item[(i)] If $m$ is odd then 

\scalebox{0.9}{
\vbox{
\[
-1 \leq \E\Big[\Big\langle\frac{X}{\Vert X \Vert_{l_2^d}},\frac{Y}{\Vert Y \Vert_{l_2^d}}
\Big\rangle^m_{l_2^d} \Big] = c_{-}(d,m)\,(1-\rho^2)^{\frac{d}{2}}\,\rho\,
{}_3 F_2\big(\frac{d+1}{2}, \frac{d+1}{2}, \frac{m+2}{2}; 
\frac{3}{2}, \frac{m+d+1}{2}; \rho^2\big) \leq 1\,,
\]
}}

where
\[
c_{-}(d,m) : = 
\frac{2}{\sqrt{\pi}}\,\frac{\Gamma\big(\frac{d+1}{2}\big)^2\,\Gamma\big(\frac{m+2}{2}\big)}
{\Gamma\big(\frac{d}{2}\big)\Gamma\big(\frac{m+d+1}{2}\big)}\,. 
\]
\item[(ii)] If $m$ is even then

\scalebox{0.9}{
\vbox{
\[
-1 \leq \E\Big[\Big\langle\frac{X}{\Vert X \Vert_{l_2^d}},\frac{Y}{\Vert Y \Vert_{l_2^d}}
\Big\rangle^m_{l_2^d} \Big] = c_{+}(d,m)\,(1-\rho^2)^{\frac{d}{2}}\,
{}_3 F_2\big(\frac{d}{2}, \frac{d}{2}, \frac{m+1}{2}; 
\frac{1}{2}, \frac{m+d}{2}; \rho^2\big) \leq 1\,,
\]
}}

where
\[
c_{+}(d,m) : = 
\frac{1}{\sqrt{\pi}}\,\frac{\Gamma\big(\frac{d}{2}\big)\,\Gamma\big(\frac{m+1}{2}\big)}
{\Gamma\big(\frac{m+d}{2}\big)}\,.
\] 
\end{itemize}
\end{thm*}
\noindent If we apply the latter result to $m = 1$ (and $d \in \N)$, we obtain a result which 
contains \cite{BOFV2014}, Lemma 2.1 as a special case if $\rho \in (-1,1)$. Clearly, that result also holds for $\rho \in \{-1,1\}$.
\begin{cor*}
Let $d \in \N$, $\rho \in [-1,1]$ and $(X^\top, Y^\top)^\top \equiv 
(X_1, \ldots, X_d, Y_1, \ldots, Y_d)^\top \sim N_{2d}\big(0, \Sigma_{2d}(\rho)\big)$. 
Then 
\[
-1 \leq \E\Big[\Big\langle\frac{X}{\Vert X \Vert_{l_2^d}},\frac{Y}{\Vert Y \Vert_{l_2^d}}
\Big\rangle_{l_2^d} \Big] = c_d\,\rho\,{}_2 F_1\big(\frac{1}{2}, \frac{1}{2},\frac{2+d}{2}; \rho^2 \big) \leq 1\,,
\]
where
\[
c_d : = \frac{1}{{}_2 F_1\big(\frac{1}{2}, \frac{1}{2}, 
\frac{2+d}{2}; 1\big)} = \frac{2}{d}\,
\frac{\Gamma\big(\frac{d+1}{2}\big)^2}{\Gamma\big(\frac{d}{2}\big)^2} 
= \frac{1}{d}\big(\E\big[\Vert X \Vert_{l_2^d}\big]\big)^2 
\]
\end{cor*} 
\begin{rem*}
This result (respectively \cite[Lemma 2.1]{BOFV2014}) shows why in the real case as 
well as the complex case the function $\text{sign}$ works so smoothly. Since if we choose 
$\text{sign}$ then we obtain in both cases an inner product (where 
$\langle x, y\rangle_{l_2^1} : = x\cdot y$ for all $x, y \in \R \equiv l_2^1$, of course). 
In particular, observe that for all $m \in \N$ $c_1 = c_{-}(1,m) = \frac{2}{\pi}$, 
$c_{+}(1,m) = 1$ and ${}_2 F_1\big(\frac{1}{2}, \frac{1}{2},\frac{1}{2}; \rho^2 \big)
= \arcsin'(\rho) = \frac{1}{\sqrt{1-\rho^2}}$ for all $\rho \in (-1,1)$. Another 
observation whch actually leads to a straightforward and short proof is the fact that in our 
calculation of the related multiple integral we also implement the probability space 
$(\Omega_d, \mathcal{B}(\Omega_d), \P_d)$, where
\[
\Omega_d : = \S^{d-1} \times \S^{d-1}, \P_d : = \underline{\sigma_d} \otimes 
\underline{\sigma_d}
\]
and $\mathcal{B}(\Omega_d)$ denotes the Borel sigma-algebra on $\Omega_d$.
\\[0.2em] 
However, in \cite{BMMN2013} the authors show that in fact 
$K_G^\R < \frac{\pi}{2\ln(1+\sqrt{2})}$, implying that in the real case $\text{sign}$ is 
not the ``optimal'' function to choose (answering a question of H. K\"{o}nig to the 
negative (see \cite{K2001})! 
\end{rem*}
\section{Completely correlation preserving functions and their impact on the upper 
bound of the Grothendieck constant}
\noindent Already while looking for the smallest upper bound of both, $K_G^\R$ and $K_G^\C$, 
we are lead to a deep interplay of different subfields of mathematics (both, pure and applied) 
including Gaussian harmonic analysis and Malliavin calculus (Mehler kernel, Ornstein-Uhlenbeck 
semigroup, Hermite polynomials, Gegenbauer polynomials\footnote{also known as \textit{ultraspherical 
polynomials}}, integration over spheres in $\R^n$), complex analysis (analytic 
continuation and biholomorphic mappings, special functions), combinatorial analysis (inversion 
of Taylor series and ordinary partial Bell polynomials), matrix analysis (positive 
semidefinite matrices, block matrices) and multivariate statistics and 
high-dimensional Gaussian dependence modelling (correlation matrices, real and complex Gaussian random vectors, Gaussian measure). 
\\[0.2em]
In particular, we have to look for those functions which map correlation matrices of any 
size and any rank entrywise into a correlation matrix of the same size again, by means of the 
so called \textit{Schur product} (also known as \textit{Hadamard product}) of matrices:
\begin{df*}[Schur product]
Let $m,n \in \N$. Let $A = (a_{ij}) \in \M(m \times n; \F)$ and $B = (b_{ij}) \in \M(m \times n; \F)$. The Schur product $A \ast B \in \M(m \times n; \F)$ is defined as
\[
(A \ast B)_{ij} : = a_{ij}\,b_{ij} \hspace{0.3cm} \big((i,j) \in [m] \times [n]\big)\,.
\] 
\end{df*}
\noindent The usefulness of the Schur product structure is reflected in the Schur product theorem which 
states that the (closed and convex) cone of all positive semidefinite matrices is stable 
under Schur multiplication:
\begin{thm*}[Schur, 1911]
Let $m,n \in \N$. Let $A = (a_{ij}) \in \M(m \times n; \F)$ and $B = (b_{ij}) \in \M(m \times n; \F)$. If both, $A$ and $B$ are positive semidefinite then also $A \ast B$ is positive semidefinite.
\end{thm*}
\noindent In particular, for any $k \in \N$ the Schur product of two $k \times k$ correlation 
matrices of the same size again is a $k \times k$ correlation matrix.  
\begin{df*}
Let $m,n \in \N$. Given $\emptyset \not= U \subseteq \F$, a function 
$f : U \longrightarrow \F$ and a matrix $A = (a_{ij}) \in \M(m \times n; U)$ put
\[
f[A] : = \big(f(a_{ij})\big) \hspace{0.3cm} \big((i,j) \in [m] \times [n]\big)\,.
\]
\end{df*}
\noindent In particular, if $f(x) = \sum_{n=0}^\infty c_n\,x^n, x \in U$, where 
$c_n \in \F$ for all $n \in \N$, we have
\[
f[A]_{ij} = \sum_{n=1}^\infty c_n\,a_{ij}^n \text{ for all } (i,j) \in [m] \times [n]\,.
\]
Functions of the latter type, where $c_n \geq 0$ for all $n \in \N$ play a significant role, 
also with respect to an analysis of $K_G^\F$. This is due to the following (cf. 
\cite[Theorem 2.1]{BGKP2016})
\begin{thm*}[Schoenberg, 1942; Rudin, 1959]
Let $f : [-1,1] \longrightarrow \R$ be a function. Then the following statements are equivalent:
\begin{itemize}
\item[(i)] $f$ is continuous and $f[A]$ is positive semidefinite for all positive semidefinite 
matrices $A$ with entries in $[-1,1]$ and of any size.
\item[(ii)] $f[A]$ is positive semidefinite for all positive semidefinite matrices $A$ with 
entries in $[-1,1]$ and of any size.
\item[(iii)] 
\[
f(x) = \sum_{n=0}^\infty a_n\,x^n \text{ for all } x \in [-1,1], \text{ where } a_n \geq 0 
\text{ for all } n \in \N_0 \text{ and } (a_n)_{n \in \N_0} \in l_1\,.
\]
\item[(iv)] $f$ is continuous from the left at $1$ and continuous from the right at $-1$ and 
can be extended to an analytic function on $\D = \{z \in \C : \vert z \vert < 1\}$ which has 
non-negative coefficients in $l_1$.      
\end{itemize}
\end{thm*}
\noindent Rudin verified that the continuity assumption in Schoenberg's proof of the 
implication ``$(i) \Rightarrow (iii)$'' can be dropped (cf. \cite{BGKP2016}, comment below 
Theorem 2.1). Their result naturally leads to the 
following definition. Given our aim to look for a common source of the real and complex case, 
we put $\F_1 : = \R$, $\F_2 : = \C$ and set $\D_d : = \{z \in \F_d : \vert z \vert < 1\} 
\hspace{0.1cm} (d \in \{1,2\})$.
\begin{df*}[Completely correlation preserving function]
Let $d \in \{1,2\}$ and $g : \overline{\D_d} \longrightarrow \F_d$ be a function.
\begin{itemize}
\item[(i)] Given $n \in \N$, $g$ is $n$-correlation-preserving (short: $n$-CP) if for any 
$n \times n$ correlation matrix $\Sigma \in C(n; \F)$ also $g[\Sigma] \in C(n; \F)$ is an 
$n \times n$ correlation matrix. 
\item[(ii)] $g$ is called completely correlation-preserving (short: CCP) if $g$ is 
$n$-correlation-preserving for all $n \in \N$.  
\end{itemize}
\end{df*}
\noindent Obviously, every CCP function $g$ satisfies $g(\overline{\D_d}) \subseteq 
\overline{\D_d}$. Since any positive semidefinite matrix factors through a correlation matrix 
(with respect to the standard matrix product), we obtain the following result, given that 
$\F_1 = \R$:
\begin{thm*}
Let $g : [-1,1] \longrightarrow \R$ be a function. Then the following statements are 
equivalent:
\begin{itemize}
\item[(i)] $g$ is CCP.
\item[(ii)] $g(1) = 1$ and $g[A]$ is positive semidefinite for all positive semidefinite matrices $A$ with entries in $[-1,1]$ and of any size.
\item[(iii)] 
\[
g(\rho) = \sum_{n=0}^\infty a_n\,\rho^n \text{ for all } \rho \in [-1,1], \text{ where } 
a_n \geq 0 \text{ for all } n \in \N_0 \text{ and } (a_n)_{n \in \N_0} \in S_{l_1}\,.
\]
\item[(iv)] $g(\rho) = \E_\P[\rho^X] = \sum_{n=0}^\infty \P\lb X = k\rb \rho^k$ is the 
probability generating function of some discrete random variable 
$X : \Omega \longrightarrow \N_0$, defined on a probability space 
$(\Omega, \mathcal{F}, \P)$.       
\end{itemize}
\end{thm*}
\noindent Built on these facts, Grothendieck's identity implies the following
\begin{cor*}[Grothendieck's identity as entrywise matrix equality]
Let $k \in \N$ and $\Sigma \in C(k; \R)$ be an arbitrarily given $k \times k$ correlation 
matrix. Then also $\frac{2}{\pi}\arcsin[\Sigma] \in C(k; \R)$ is a 
$k \times k$ correlation matrix. There exist column vectors $u_1, u_2, \ldots, u_k \in 
\S^{k-1}$ such that for any $X \equiv (X_1, \ldots, X_k)^\top \sim N_k(0,I_k)$ 
\[
\frac{2}{\pi}\arcsin[\Sigma] = \E\big[\Theta(u)\big]\,,
\]
where the random matrix $\Theta(u) \equiv \Theta(vec(u_1, \ldots, u_k))$ is given as
\[
\Theta(u)_{ij} : = {\text{sign}}(\langle X, u_i\rangle){\text{sign}}(\langle X, u_j\rangle) 
\]
for all $i, j \in [k]$. $\Theta(u) = {\text{sign}}\Big[\big(\bigoplus_{i=1}^k X^\top\big)u\Big]
{\text{sign}}\Big[\big(\bigoplus_{i=1}^k X^\top\big)u\Big]^\top$ is a random correlation matrix of rank $1$.
\end{cor*}
\noindent Let $d \in \{1,2\}$ and put $\T_d : = \partial\,\D_d : = 
\{z \in \F_d : \vert z \vert = 1\}$. Another important estimation (even with upper bound $1$) 
which might support our search for a ``suitable'' CCP function which is different from the CCP 
function $\frac{2}{\pi}\arcsin$ is the following one  
\begin{prp*}
Let $d \in \{1,2\}$, $z \in \overline{\D_d}$ and $b : \F_d \longrightarrow \T_d$. Let 
$h_b : \overline{\D_d} \longrightarrow \overline{\D_d}$ be a function such that
\[
h_b(z) = \E\big[b(Z)\,\overline{b(W)}\big]
\]
for any $z \in \overline{\D_d}$ and any $(Z, W)^\top \sim 
{\F_d}N_{2}\big(0, \Sigma_{2}(z)\big)$. Then
\begin{itemize}
\item[(i)] For all $k \in \N$ and $\Sigma \in C(k; \F_d)$, 
\[
h_b[\Sigma] \in C(k; \F_d)\,.
\]
\item[(ii)] For all $k \in \N$, $\Sigma \in C(k; \F_d)$ and $A \in \M(k \times k; \F_d)$,
\[
\big\vert{\text{tr}}\big(A^\ast\,{}{h_b\big[\Sigma\big]}\big)\big\vert \leq
\max_{\Theta \in C_1(k; \F_d)}\big\vert{\text{tr}}\big(A^\ast\,\Theta\big)\big\vert\,.
\]
\end{itemize}
Moreover, $h_b(1)=1$ and $h_b(0) = \big\vert\E[b(Z)]\big\vert^2$ for any $Z \in 
{\F_d}N_{1}\big(0, 1\big)$.
\end{prp*}
\begin{rem*}
In the statistical learning community in artificial intelligence a function $b : \R^k 
\longrightarrow \{-1,1\}, k \in \N$ is called ``concept'' (cf. e.g. \cite{S2009}). We adopt 
this name. If the function $h_b$ were invertible, for some $b$ then $h_b^{-1}$ cannot 
be a CCP function (else $K_G^\F \leq 1$\,(!) - a contradiction). In particular, $h_b^{-1}$ 
cannot be represented as a power series with non-negative coefficients, such as e.\,g. 
$h_{\text{sign}}^{-1}(y) = 
\sin\big(\frac{\pi}{2}\,y\big) = \sum_{n=0}^\infty (-1)^n 
\frac{\pi^{2n+1}}{2^{2n+1}\,(2n+1)!}\,y^{2n+1}$, $y \in [-1,1]$ if we consider
$h_{\text{sign}}(x) = \frac{2}{\pi}\arcsin(x)$, $x \in [-1,1]$.
\end{rem*}
\noindent However, we will see that inversion of CCP functions plays the key role regarding 
the search for the lowest upper bound of the Grothendieck constant $K_G^\F$. Unfortunately, a 
closed form representation of the coefficients of the inverse of a Taylor series runs against 
a well-known combinatorial complexity issue (due to the presence of ordinary partial 
Bell polynomials as building blocks of these coefficients - cf. research problem 7.1 below 
for details), which in 
general does not allow a closed form representation of these coefficients, such as is the 
case with the inverse of Haagerup's function $\overline{\D} \ni z \mapsto 
\frac{\pi}{4}\,z\,{}_2 F_1\lb\frac{1}{2}, \frac{1}{2}, 2; \vert z \vert^2 \rb$ in the complex 
case (cf. \cite{H1987}, Remark on page 216), as opposed to Grothendieck's function 
$[-1,1] \ni \rho \mapsto \frac{2}{\pi}\,\rho\,{}_2 F_1\lb\frac{1}{2},\frac{1}{2}, \frac{3}{2}; 
\rho^2 \rb = \frac{2}{\pi}\arcsin(\rho)$ in the real case. Clearly, further research is here 
required (as indicated below). 
\\[0.2em]
Keeping inversion of CCP functions in mind we now are going to present two key results of our 
research.
\begin{thm*}
Let $d \in \{1,2\}$ and $m, n \in \N$. Let $A \in \M(m \times m; \F_d)$ and 
$B \in \M(n \times n; \F_d)$ be positive semidefinite and $Z \in \M(m \times n; \F_d)$.
Let $f : \overline{\D_d} \longrightarrow \F_d$ and $g : \overline{\D_d} \longrightarrow \F_d$ 
be functions. Suppose that $f$ can be represented as $f(z) = \sum_{n=1}^\infty a_n \,z^n$ for 
all $z \in \overline{\D_d}$, where $a_n \geq 0$ for all $n \in \N$. Assume that $g$ can be 
represented as $g_b(z) = \sum_{n=1}^\infty b_n \,z^n$ for all $z \in \overline{\D_d}$, where 
$b_n \in \F_d$ for all $n \in \N$. Assume further that
\[
\vert b_n \vert \leq a_n \text{ for all } n \in \N\,.
\]  
If the block matrix 
\[
\Sigma : =
\begin{pmatrix}
A & Z \\
Z^\ast & B
\end{pmatrix}
\] 
is positive semidefinite, and if all matrices $A, B$ and $Z$ have entries in $\overline{\D_d}$ 
then also the block matrix
\[
\Sigma_{f,g} : = 
\begin{pmatrix}
f[A] & g[Z] \\
g[Z]^\ast & f[B]
\end{pmatrix}
\]
is positive semidefinite. If in addition $f(r) = 1$ for some $0 < r \leq 1$, and if 
$\Sigma \in C(m+n, \F_d)$ is a correlation matrix then also 
$(r\Sigma)_{f, g} \in C(m+n, \F_d)$ is a correlation matrix.      
\end{thm*}
\noindent To prepare the underlying ideas of our next result, we carefully list the single 
steps and assumptions, possibly leading to an algorithmic approach regarding the 
implementation of an approximation to the lowest upper bound of the Grothendieck constant 
$K_G^\F$. All of these conditions are satisfied for the Grothendieck function and the Haagerup 
function and were used to construct a ``small'' upper bound of the respective 
Grothendieck constant. However, if we wish to apply these steps to a ``better 
fitting function'' $b$ which is different from both, the Grothendieck function and the 
Haagerup function, we are strongly confronted with non-trivial persistent combinatorial issues; described in Section 5 below.
\\[0.2em]
So, fix $d \in \{1,2\}$ and consider the following workflow step by step.
\begin{center}
\setlength{\fboxsep}{25pt}
\fbox{\parbox{0.99\columnwidth}{
\begin{itemize}
\item[(SIGN)] Choose a function $b : \F_d^k \longrightarrow \T_d$ for some ``suitable'' 
$k \in \N$ and consider its allocated CCP function $h_b$, constructed according to the lines 
of the above Proposition (such as e.g. the Grothendieck function, respectively the Haagerup 
function $h_{\text{sign}}$). 
\item[(H)] Assume that $h_b(0) = 0$ and that
$h_b : \overline{\D_d} \longrightarrow \overline{\D_d}$ is a homeomorphism.
\item[(RA)] 
Assume that $g_b : = h_b^{-1} : \overline{\D_d} \longrightarrow \overline{\D_d}$ can be 
represented as $g_b(z) = \sum_{n=1}^\infty \beta_n \,z^n$ for all $z \in \overline{\D_d}$, 
where $\beta_n \in \R$ (!) for all $n \in \N$ and $\beta \equiv (\beta_n)_{n \in \N} 
\in l_1$; i.e., $\Vert \beta \Vert_1 \equiv \sum_{n=1}^\infty \vert \beta_n \vert < \infty$.
\item[(ABS)] Put $f_b(z) : = \sum_{n=1}^\infty \vert \beta_n \vert \,z^n \hspace{0.1cm} 
(z \in \overline{\D_d})$.
\item[(PI(1))] Assume that $f_b(r) = 1$ for some $0 < r < 1$.
\end{itemize}
}}
\end{center}
Clearly, we have $1 = g_b(1) = \sum_{n=1}^\infty \beta_n \leq \Vert \beta \Vert_1 = f_b(1)$.
Moreover, $\vert f_b(rz) \vert \leq f_b(r) = 1$ for any $z \in \overline{\D_d}$, implying 
that $\overline{\D_d} \ni h_b\circ f_b(rz) = h_b(f_b(rz))$ is well-defined if $z \in 
\overline{\D_d}$.   
Note that $r$ depends on the choice of $b$. 
\\[0.2em]
Let $k \in \N$. Fix an arbitrary correlation matrix 
\[
\Sigma = \begin{pmatrix}
A & Z \\
Z^\ast & B
\end{pmatrix} \in C(2k, \F_d)\,.
\]
\noindent Combining the previous two results and since $g_b(\overline{z}) = \overline{g_b(z)}$ 
it therefore follows that
\[
\Sigma_b(r) : = \begin{pmatrix}
h_b\circ f_b[r\,A] & r\,Z \\
r\,Z^\ast & h_b\circ f_b[r\,B]
\end{pmatrix} = h_b \lsb 
\begin{pmatrix}
f_b[r\,A] & g_b[r\,Z] \\
(g_b[r\,Z])^\ast & f_b[r\,B]
\end{pmatrix} 
\rsb \in C(2k, \F_d)
\]  
again is a correlation matrix. Consequently, we have
\[
r\,\vert{\text{tr}}\big(J(A)\,\Sigma)\vert = \vert{\text{tr}}\big(J(A)\,(r\,\Sigma)\big)\vert 
= \vert{\text{tr}}\big(J(A)\,\Sigma_b(r)\big)\vert \leq K_G^\F\, 
\Vert A \Vert_{\infty, 1}\,,
\]
and hence   
\begin{thm*}
Let $k \in \N$, $d \in \{1,2\}$, $b : \F_d^k \longrightarrow \T_d$ and $h_b : \overline{\D_d} 
\longrightarrow \overline{\D_d}$ the allocated CCP function. Assume that the assumptions (H) 
and (RA) of the workflow hold. Let $f_b : \overline{\D_d} \longrightarrow \F_d$ be 
constructed as above and assume that $f_b$ satisfies the condition (PI(1)) for some $0 < r^
\ast < 1$. Then 
\[
K_G^\F \leq \frac{1}{r^\ast}\,.
\]
\end{thm*}
\begin{ex*}[Krivine's upper bound reproduced]
Let $\F = \R$. Consider $b : = \text{sign}$. Due to Grothendieck's identity we know that
\[
h_b(\rho) = \frac{2}{\pi}\arcsin(\rho) \text{ for all } \rho \in [-1,1]\,.
\]
The continuous function $h_b : [-1,1] \longrightarrow [-1,1]$ is strictly increasing and hence 
invertible, with continuous inverse 
\[
g_b(\tau) : = h_b^{-1}(\tau) = \sin\big(\frac{\pi}{2}\,\tau\big) = \sum_{n=0}^\infty (-1)^n\,
\frac{\pi^{2n+1}}{2^{2n+1}\,(2n+1)!}\,\tau^{2n+1} \hspace{0.2cm} (\tau \in [-1,1])\,.
\] 
Hence,
\[
f_b(\tau) = \sum_{n=0}^\infty \frac{\pi^{2n+1}}{2^{2n+1}\,(2n+1)!}\,\tau^{2n+1} = 
\frac{1}{i}\sin\big(\frac{\pi}{2}\,i\,\tau\big) = \sinh\big(\frac{\pi}{2}\,\tau\big)  
\hspace{0.2cm} (\tau \in [-1,1])\,.
\]
Since $f_b\lb\frac{2\ln(1+\sqrt{2})}{\pi}\rb = 1$ it follows that
\[
K_G^\R \leq \frac{\pi}{2\ln(1+\sqrt{2})}\,. 
\]
\end{ex*}
\section{Emerging research problems}
\noindent Not very surprisingly, the long-standing, intensive and technically quite demanding 
attempts to detect the - still not available - value of the both Grothendieck constants (open 
since 1953) leads to further challenging tasks and open problems, such as the following ones; 
addressed in particular to highly motivated students who also wish to get a better 
understanding of the reasons underlying these difficulties.   
\subsection{Research problem 1: Grothendieck's constant versus Taylor series inversion}
\noindent Only between 2011 and 2013 it was shown (cf. \cite{BMMN2013}) 
that $K_G^\R$ is strictly smaller than Krivine's upper bound, stating that 
$K_G^\R < \frac{\pi}{2\ln(1+\sqrt{2})}$. 
Consequently, in the real case $\text{sign}$ is not the ``optimal'' function to choose 
(answering a question of H. K\"{o}nig to the negative - cf. \cite{K2001}). So, if we wish to 
reduce the value of the upper bound of the real Grothendieck constant we have to look for 
functions $b : \R^k \longrightarrow \{-1,1\}$ which are different from 
$\text{sign} : \R \longrightarrow \{-1,1\}$. 
However, these functions should satisfy all of the conditions in the listed workflow above. 
In particular, we have to look for both, the coefficients of the Taylor series of $h_b$ and 
the coefficients of the Taylor series of the \textit{inverse} function $g_b : = h_b^{-1}$. It 
is well-known that the latter task increases strongly in computational complexity if we want 
to calculate such Taylor coefficients of a higher degree, leading to the involvement of non-
trivial combinatorial facts, reflected in the use of partitions of positive integers and 
partial exponential Bell polynomials as part of the Taylor coefficients of the inverse Taylor 
series (a thorough introduction to this framework including the related 
Lagrange-B\"urmann inversion formula is given in \cite{K1951, C1974}).
\\[0.2em]
To reveal the origin of these difficulties let us focus on the real case, with $k=1$. First note that 
\[
h_b(\rho) = \sum_{n=1}^\infty \langle b, H_n \rangle_{\gamma_1}^2\,\rho^n
\]
for all $\rho \in [-1,1]$, where for $n \in \N$ and $x \in \R$
\[
H_n(x) : = \frac{1}{\sqrt{n!}} (-1)^n \exp\lb \frac{x^2}{2}\rb\,\frac{d^n}{dx^n}
\exp\lb-\frac{x^2}{2} \rb 
\]
denotes the (probabilistic version of the) $n$-th Hermite polynomial and $\gamma_1$ the 
Gaussian measure on $\R$ (cf. e.g. \cite{B1998}). Put $0 \leq \alpha_n : =
\langle b, H_n \rangle_{\gamma_1}^2$. If $h_b'(0) \not= 0$ we know that at least the 
real-analytic function $h_b \mid (-1,1)$ is invertible around $0 = h_b(0)$. Its inverse is 
also expressible as a power series there; i.e., around $0$, $\big(h_b \mid (-1,1)\big)^{-1}$ 
is real-analytic, too. Hence, given the assumption (RA) it follows that $h_b^{-1}(y) = 
g_b(y) = \sum_{n=1}^\infty \beta_n \,y^n$ for all $y \in [-1,1]$, where $\beta_1 = 
\frac{1}{\alpha_1}$ and
\begin{eqnarray*}
\beta_n & = & \frac{1}{n}\,\sum_{k=1}^{n-1}\frac{1}{\alpha_1^{n+k}}(-1)^{k}\,
\binom{n-1+k}{k}\,B_{n-1,k}^{\circ}\big(\alpha_2, \alpha_3, \ldots, \alpha_{n-k+1}\big)\\
& = & \frac{1}{n}\,\sum_{k=1}^{n-1}(-1)^{k}\,\binom{n-1+k}{k}\,
B_{n-1,k}^{\circ}\big(\alpha_2, \alpha_3\alpha_1, \alpha_4\alpha_1^2 \ldots, 
\alpha_{n-k+1}\alpha_1^{(n-1)-k}\big)
\end{eqnarray*}
for all $n \in \N_2$. Thereby,
\[
B_{n,k}^{\circ}\big(x_1, x_2, \ldots, x_{n+1-k}) : =
\sum_{\nu \in P(n, k)}\frac{k!}{\prod_{i=1}^{n+1-k}{\nu_i}!}\,\prod_{i=1}^{n+1-k} 
x_i^{\nu_i} = \sum_{\nu \in P(n, k)}k!\frac{x^\nu}{\nu!}
\] 
denotes the ordinary partial Bell polynomial and $P(n, k)$ indicates the set of all 
multi-indices $\nu \equiv (\nu_1, \nu_2, \ldots, \nu_{n+1-k}) \in \N_0^{n+1-k}$ ($k \leq n$) 
which 
satisfy the Diophantine equations $\sum_{i=1}^{n+1-k} \nu_i = k$ and $\sum_{i=1}^{n+1-k} 
i\,\nu_i = n$; i.e., summation is extended over all partitions of the number $n$ into 
exactly $k$ summands (cf. e.g. \cite{C1974, C2011, MMKS2019}). We explicitly list $\beta_2, \beta_3, \beta_4, \beta_5, \beta_6$ 
and $\beta_7$:
\[
\beta_2\alpha_1^3 = -\alpha_2\,,
\]
\[
\beta_3\alpha_1^5 = -\alpha_1\alpha_3 + 2\alpha_2^2\,,
\]
\[
\beta_4\alpha_1^7 = -\alpha_1^2\alpha_4 + 5\,\alpha_1\alpha_2\,\alpha_3 - 5\,\alpha_2^3\,,
\]
\[
\beta_5\alpha_1^9 = -\alpha_1^3\alpha_5 + 6\,\alpha_1^2\alpha_2\,\alpha_4 + 
3\,\alpha_1^2\alpha_3^2 - 21\,\alpha_1\alpha_2^2\,\alpha_3 + 14\,\alpha_2^4\,,
\]
\[
\beta_6\alpha_1^{11} = -\alpha_1^4\alpha_6 + 7\,\alpha_1^3\alpha_2\,\alpha_5 + 
7\,\alpha_1^3\alpha_3\,\alpha_4 - 28\,\alpha_1^2\alpha_2\,\alpha_3^2 - 
28\,\alpha_1^2\alpha_2^2\,\alpha_4 + 84\,\alpha_1\alpha_2^3\,\alpha_3 - 42\,\alpha_2^5\,,
\]
\scalebox{0.9}{
\vbox{
\begin{eqnarray*}
\beta_7\alpha_1^{13} & = & -\alpha_1^5\alpha_7 + 8\,\alpha_1^4\alpha_2\,
\alpha_6 + 8\,\alpha_1^4\alpha_3\,\alpha_5 + 4\,\alpha_1^4\alpha_4^2 - 
36\,\alpha_1^3\alpha_2^2\,\alpha_5 - 72\,\alpha_1^3\alpha_2\,\alpha_3\,\alpha_4\\
& - & 12\,\alpha_1^3\alpha_3^3 + 120\,\alpha_1^2\alpha_2^3\,\alpha_4 + 180\,\alpha_1^2\alpha_2^2\,\alpha_3^2 - 330\,\alpha_1\alpha_2^4\,\alpha_3 + 132\,\alpha_2^6\,.
\end{eqnarray*}
}}

\noindent It appears to us that there is a general pattern in these formulas which might lead to an expression of the following type:

\scalebox{0.82}{
\vbox{
\begin{equation*}
n\beta_n\,\alpha_1^{2n-1} = - n\alpha_1^{n-2}\alpha_n + 
\sum_{k=2}^{n-2}(-1)^k\,m_k\,\binom{(n-1)+k}{k}\,\alpha_1^{(n-1)-k}\,
p_k(\alpha_2, \alpha_3, \ldots, \alpha_n) + (-1)^{n-1}\binom{2(n-1)}{n-1}\alpha_2^{n-1}\,, 
\end{equation*}
}}

\noindent where $m_k \in [k]$, $p_k(\alpha_2, \alpha_3, \ldots, \alpha_n) : = 
\sum\limits_{(\nu_2,\ldots,\nu_n)\in \A_k}{\prod_{l=2}^n}\alpha_l^{\nu_l}$ and 
$\A_k \subset \big\{\mu \in \N_0^{n-1} : \sum_{l=2}^n l\,\mu_{l-1} = n-1+k\big\}$.
\\[0.2em]
The strong difficulties are twofold: already in the one-dimensional case we need to know 
the explicit value of all of the Fourier coefficients 
\[
\langle b, H_n \rangle_{\gamma_1} = 
\frac{1}{\sqrt{2\pi}}\int_\R b(x) H_n(x)\,\exp\big(-\frac{1}{2}\,x^2\big)\,dx = 
\E[b(X)H_n(X)] = \frac{d^n}{dt^n}\E[b(X+t)]\Big\vert_{t=0}\,,
\]
where $X \sim N_1(0,1)$, together with a closed form expression (if available at all) of all 
of the coefficients $\beta_n$, where the latter involves (in general, unknown values of) 
ordinary partial Bell polynomials. It appears to us that in general one cannot use proofs by 
standard induction on $n \in \N$ to verify statements about Bell polynomials. Here, the 
Noetherian Induction Principle seems to be more appropriate (cf \cite{OeBELL2020}). 
\\[0.2em]
At least in the case of Grothendieck's $b : = \text{sign} = 2 \ind_{[0, \infty)} - 1$ we can show that for all $n \in \N$, 
\[
\alpha_{2(n-1)} = 0, \alpha_1 = \frac{2}{\pi} \text{ and } \alpha_{2n+1} = 
\frac{2}{\pi}\,\frac{((2n-1)!!)^2}{(2n+1)!}
\]
Consequently, if $n \in \N$, the (best-known) power series representation of 
$g_b(y) = \sin\lb\frac{\pi}{2}\,y\rb$ leads to the following interesting identity:

\scalebox{0.85}{
\vbox{
\begin{eqnarray*}
(-1)^{n}\,\frac{\big(\pi/2\big)^{2n+1}}{(2n+1)!} & = & \beta_{2n+1}\\
& = & \frac{\big(\pi/2\big)^{2n+1}}{(2n+1)!}\,\sum_{k=1}^{2n}(-1)^{k}\,\frac{(2n+k)!}{k!}\,
B_{2n,k}^{\circ}\lb 0, \frac{1}{6}, 0, \frac{3}{40}, 0, \frac{5}{112}, \ldots, 
\frac{1+(-1)^{k+1}}{2}\frac{((2n-k)!!)^2}{(2n-k+2)!}\rb\,.
\end{eqnarray*}
}}

\noindent However, already in this case we don't know a closed form expression for the numbers 
$B_{2n,k}^{\circ}\lb 0, \frac{1}{6}, 0, \frac{3}{40}, 0, \frac{5}{112}, \ldots, 
\frac{1+(-1)^{k+1}}{2}\frac{((2n-k)!!)^2}{(2n-k+2)!}\rb$. An even stronger problem appears in 
the complex case, since already a closed-form formula for the coefficients of the Taylor 
series of the inverse of the Haagerup function is still unknown (cf. \cite{H1987}). Here, we 
would like to list the very 
recent paper \cite{QNLY2020}, where the authors point to similar difficulties including the 
formulation of related - open - problems. Moreover, the solved examples in \cite{QNLY2020} 
show the large barriers which we have to resolve while working with (partial) Bell 
polynomials. Here, \cite{OeBELL2020} will 
disclose further surprising properties of these polynomials. In particular, we will provide a 
closed form sum representation of the polynomials underlying a 
very useful - recursive - construction of T. M. Apostol (see \cite{A2000}) of the 
coefficients of the Taylor series of inverse functions.   
\\[0.2em]
Therefore, given the intrinsic combinatoric and computational complexity regarding the 
determination of the Taylor series coefficients of the Taylor series of the non-CCP function 
$g_b$ (via the Fa\'a-di Bruno formula), related research problems (which actually do not 
require any knowledge of the Grothendieck inequality) could be the following ones:
\begin{itemize}
\item[(1-RP1)] Let $n \in \N_4$ and $k \in \{2,3, \ldots, n-2\}$. Recall that

\scalebox{0.82}{
\vbox{
\begin{equation*}
n\beta_n\,\alpha_1^{2n-1} = - n\alpha_1^{n-2}\alpha_n + 
\sum_{k=2}^{n-2}(-1)^k\,m_k\,\binom{(n-1)+k}{k}\,\alpha_1^{(n-1)-k}\,
p_k(\alpha_2, \alpha_3, \ldots, \alpha_n) + (-1)^{n-1}\binom{2(n-1)}{n-1}\alpha_2^{n-1}\,. 
\end{equation*}
}}

Prove this representation and determine $\A_k$ and the numbers $m_k$ therein explicitly 
(if feasible at all)!
\item[(1-RP2)] Continue to investigate the structure of partial Bell polynomials; 
possibly with the aid of high-performance computers and related (algebraic) software tools.
\end{itemize}
\subsection{Research problem 2: Grothendieck's inequality and copulas}
\noindent If we thouroughly overhaul the CCP function $[-1,1] \ni \rho \mapsto 
\frac{2}{\pi}\arcsin(\rho)$ we recognise that some knowledge of Gaussian copulas (i.e., 
finite-dimensional multivariate distribution functions of univariate marginals generated by 
the distribution function of Gaussian random vectors - cf. e.g. \cite{Sk1959, 
N2006, UADF2017} and \cite{Oe2015}) and (the probabilistic version) of the Hermite 
polynomials might become very fruitful regarding our indicated search for different ``
suitable'' CCP functions. $[-1,1] \ni \rho \mapsto \psi(\frac{1}{2},\frac{1}{2}; t) = 
\frac{2}{\pi}\arcsin(\rho)$ namely reveals as a special case of the CCP function

\scalebox{0.85}{
\vbox{ 
\[
[-1,1] \ni \rho \mapsto \psi(p,p; \rho) = \frac{1}{c(p)}\sum_{n=1}^\infty \frac{1}{n}\,
H_{n-1}^2(\Phi^{-1}(p))\,\rho^{n} = \frac{1}{2\pi\,p(1-p)}\,\exp\big(-(\Phi^{-1}(p))^2\big)\,
\rho\sum_{n=0}^\infty \frac{1}{n+1}\,H_n^2(\Phi^{-1}(p))\,\rho^{n}\,, 
\]
}}

where $0 < p < 1$ and
\[
c(p) : = \frac{p(1-p)}{\varphi^2(\Phi^{-1}(p))} = 2\pi\,p(1-p)\exp\big((\Phi^{-1}(p))^2\big) 
= \sum_{n=0}^\infty \frac{1}{n+1}\,H_{n}^2(\Phi^{-1}(p))\,.
\]
If we put $b_p : = 2\ind_{[\Phi^{-1}(p), \infty)} -1 = 1 - 2\ind_{(-\infty, \Phi^{-1}(p))} 
\in \{-1,1\}$, then the tetrachoric series expansion of the bivariate Gaussian copula (cf. 
\cite{G1963, M2013, AF2014}) implies the following generalisation of Grothendieck's 
identity:
\begin{eqnarray*}
h_p(\rho) \, : = \, h_{b_p}(\rho) : & = & \E\big[b_p(X)\,b_p(Y)\big]
\, = \, (2p-1)^2 + \frac{2}{\pi}\,\exp\big(-(\Phi^{-1}(p))^2\big)\,
\sum_{n=1}^\infty \frac{1}{n}\,H_{n-1}^2(\Phi^{-1}(p))\rho^{n}\\
& = & (2p-1)^2 + 4p(1-p)\psi(p,p; \rho)\,.
\end{eqnarray*}
Due to our construction of $\psi(p,p; \cdot)$ the latter is clearly equivalent to
\[
\rho(b_p(X),b_p(Y)) = \psi(p,p; \rho)
\]
for all $p \in (0,1)$, $\rho \in [-1,1]$ and $(X,Y) \sim N_2(0, \Sigma_2(\rho))$, where
$\rho(b_p(X),b_p(Y))$ denotes Pearson's correlation coefficient between the random variables
$b_p(X)$ and $b_p(Y)$. Unfortunately,
\[
h_p(\rho) = \psi(p,p; \rho) \text{ for all } \rho \in [-1,1] \text{ iff } p = \frac{1}{2}\,.
\]
\\[0.2em]
These facts imply the following research problems at least:
\begin{itemize}
\item[(2-RP1)] Prove whether there are $p \in (-1,1)\setminus\{\frac{1}{2}\}$ and functions 
$\chi_p : \R \longrightarrow \{-1,1\}$ such that $\psi(p,p; \rho) = 
h_{\chi_p}(\rho) = \E\big[\chi_p(X)\,\chi_p(Y)\big]$ for all $\rho \in [-1,1]$ 
and $(X,Y) \sim N_2(0, \Sigma_2(\rho))$, so that the condition (SIGN) of our workflow is 
satisfied for $h_{\chi_p}$.
\item[(2-RP2)] Generalise the above approach (which is built on the tetrachoric series of the 
bivariate Gaussian copula) to the $n$-variate case, where $n \in \N_3$ and adapt problem 
(2-RP1) accordingly.
\item[(2-RP3)] Verify whether the above approach can be transferred to the complex case. Could 
we then similarly generalise the Haagerup identity?
\item[(2-RP4)] If (2-RP1), respectively (2-RP2) holds, prove whether the remaining 
conditions (H), (RA), (ABS) and (PI1) of the workflow hold. If this were the case 
calculate (respectively approximate numerically) the related upper bound of $K_G^\R$. Include 
high-performance computers and computer algebra systems if necessary. 
\end{itemize}
\subsection{Research problem 3: Grothendieck's inequality and\\ non-commutative dependence 
structures in quantum mechanics}
\noindent Even a mathematical modelling of \textit{non-commutative dependence} 
in quantum theory and its applications to quantum information and quantum computation is strongly linked with the existence of the real Grothendieck constant $K_G^\R$. 
\\[0.2em]
The latter can be very roughly adumbrated as follows: the experimentally proven non-
Kolmogorovian (non-commutative) nature of the underlying probability theory of quantum 
physics leads to the well-known fact that in general a normal state of a composite quantum 
system cannot be represented as a convex combination of a product of normal states of the 
subsystems. This phenomenon is known as \textit{entanglement} or \textit{quantum 
correlation}. The Einstein-Podolsky-Rosen paradox, the violation of Bell's inequalities 
(limiting spatial correlation) and the Leggett-Garg inequalities (limiting temporal 
correlation) in quantum mechanics and related theoretical and experimental research implied a 
particular focus on a deeper understanding of this type of correlation - and hence to the 
\textit{modelling of a specific type of dependence} of two (ore more) quantum observables in 
a composite quantum system, measured by two (or more) space-like separated instruments, each 
one having a classical parameter (such as the orientation of an instrument which measures the 
spin of a particle). The transition probability function, i.\,e., the \textit{joint} 
probability distribution of observables in some fixed state of the system (considered as a 
function of the above-mentioned parameters) may violate Bell's inequalities and is therefore 
not realisable in ``classical'' (commutative) physics. The surprising fact, firstly 
recognised by B.S. Tsirelson (cf. \cite{T1980, T1993} and \cite{P2018}), is that also this - 
experimentally verified - gap is an implication of the existence of the real Grothendieck 
constant $K_G^\R > 1$\,!\footnote{also known as \textit{Tsirelson bound}} In other words, 
$K_G^\R$ indicates ``how non-local quantum mechanics can be at most''.    
\\[0.2em]
Already in the classical Kolmogorovian model, i.\,e., in the framework of probability space 
triples $(\Omega, \mathcal{F}, \P)$, a rigorous description of tail dependence - which \textit
{exceeds} the standard dependence measure, given by Pearson's correlation coefficient, is a 
challenging task. To disclose (and simulate) the geometry of dependence one has to 
determine finite-dimensional multivariate distribution functions of univariate marginals, 
hence \textit{copulas}. In the description of research problem 2 we have seen that Gaussian 
copulas are lurking in Grothendieck's identity. More precisely, we have (cf. \cite{S1889}):  
\begin{ex*}[Stieltjes, 1889]
Let $\rho \in [-1,1]$. Let $X, Y \sim N_1(0,1)$ such that $\E[XY] = \rho$. Then
\[
\E[{\text{sign}}(X){\text{sign}}(Y)\big] = {4\,C\big(\frac{1}{2}, \frac{1}{2}; \rho \big) - 1} = {\frac{2}{\pi}\arcsin(\rho)} = \frac{2}{\pi}\arcsin\lb\E[XY]\rb,
\]
where $[-1,1] \ni \rho \mapsto C\lb\frac{1}{2}, \frac{1}{2}; \rho \rb$ 
denotes the bivariate Gaussian copula with Pearson's correlation coefficient $\rho$ as parameter, evaluated at $(\frac{1}{2}, \frac{1}{2})$. 
\end{ex*}
\noindent Keeping a \textit{non-commutative} version of Grothendieck's inequality at the back 
of our mind (cf. \cite{P2012, T1980, T1993}) these facts might lead to problems of the 
following type:
\begin{itemize}
\item[(3-RP1)] Look for objects like ``non-commutative copulas'', leading to a search for ``non-
commutative distribution functions''.
\item[(3-RP2)] Create a ``multivariate'' spectral theory of \textit{non-commuting} normal 
operator tuples and introduce non-commutative tail dependency measures in non-commutative 
$C^\ast$-algebras and operator spaces.
\end{itemize}
\noindent We finish this research list with a few completely open questions and a conclusion:
\\[0.2em]
Can we remove the underlying Gaussian structure in Grothendieck's inequality (for both 
fields, $\R$ and $\C$) and implement tail dependent distribution functions instead (such as 
the generalized extreme value (GEV) distribution) \text{and} maintain the inequality? If this 
were the case, could that approach also be used to improve the lower and the upper bound of 
the Grothendieck constant? What about infinitely divisible probability distributions in 
general? It is very likely that just the use of correlation matrices in the trace inequality version of GT would then no longer suffice.
\subsection{Conclusion}
\noindent Apart from its comprehensive mathematical fascination and richness the highly fascinating open problem of determining the value of the both Grothendieck constants 
$K_G^\R$ and $K_G^\C$ meanwhile is deeply linked with different fields, even including applications in computer science, high-dimensional data analyis and quantum mechanics. Therefore, related research very likely would lead to a very fruitful exchange of related information and collaboration with experts from different fields, working in theory and in practice. 
\\[0.2em]
\madd{\textbf{Any comments, remarks, questions, ideas, corrections and suggestions are highly 
welcome!}}
\begin{ack*} I would like to thank Feng Qi for very helpful comments including the correction 
of a hidden typo in Section 5.1. 
\end{ack*}


\end{document}